\documentclass[reqno]{amsart}

\usepackage{amsmath}
\usepackage{amscd}
\usepackage{amssymb}
\usepackage{enumerate}
\usepackage{amsfonts}
\usepackage{mathrsfs}

\DeclareFontEncoding{OT2}{}{}

\numberwithin{equation}{section}

\newcommand{\im}{\operatorname{Im}}
\newcommand{\re}{\operatorname{Re}}

\newcommand{\rar}[1]{\stackrel{#1}{\longrightarrow}}

  \newcommand{\ga}{\gamma}
\newcommand{\Ga}{\Gamma} \newcommand{\de}{\delta} 
\newcommand{\la}{\lambda} 
 
\newcommand{\eps}{\epsilon} \newcommand{\sg}{\sigma}
  
\newcommand{\om}{\omega}

\newcommand{\bC}{{\mathbb C}}

\newcommand{\bR}{{\mathbb R}}

\newcommand{\cE}{{\mathcal E}}

\newcommand{\cJ}{{\mathcal J}}

\newcommand{\cL}{{\mathcal L}}

\newcommand{\cN}{{\mathcal N}}

\newcommand{\cS}{{\mathcal S}}

\newcommand{\cV}{{\mathcal V}}
\newcommand{\cW}{{\mathcal W}}

\newcommand{\sN}{{\mathscr N}}

\newcommand{\sS}{{\mathscr S}}

\newcommand{\fB}{{\mathfrak B}}

\newcommand{\fg}{{\mathfrak g}}

\newcommand{\fm}{{\mathfrak m}}
\newcommand{\fn}{{\mathfrak n}}

\newcommand{\fs}{{\mathfrak s}}

\newcommand{\Eb}{\overline{E}}

\newcommand{\Lb}{\overline{L}}

\newcommand{\Xb}{\overline{X}}
\newcommand{\Yb}{\overline{Y}}

\newcommand{\dd}{\partial}

\newcommand{\barr}[1]{\overline{#1}}

\newcommand{\Aut}{\operatorname{Aut}}

\newcommand{\pr}{\mathrm{pr}}

\newcommand{\rk}{\operatorname{rk}}

\newcommand{\Wedge}{\bigwedge}

\newcommand{\sbr}{\smallbreak}
\newcommand{\mbr}{\medbreak}

\newcommand{\Ann}{\operatorname{Ann}}

\newcommand{\cou}[1]{\bigl[#1\bigr]_{\mathrm{cou}}}

\newcommand{\bra}{\bigl\langle}
\newcommand{\ket}{\bigr\rangle}

\newcommand{\pair}[1]{\bra #1\ket}

\newcommand{\res}[1]{\big\vert_{#1}}

\newcommand{\matr}[4]{\left(\begin{array}{cc} \! #1 & \! #2 \! \\ \! #3 & \! #4 \! \end{array}\right)}

\newtheorem{ithm}{Theorem}

\newtheorem{thm}{Theorem}[section]

\newtheorem{lem}[thm]{Lemma}
\newtheorem{prop}[thm]{Proposition}

\theoremstyle{remark}
\newtheorem{rem}[thm]{Remark}

\newcommand{\tc}[1]{T_{\bC}#1}
\newcommand{\tcs}[1]{T_{\bC}^*#1}
\newcommand{\tts}[1]{T#1\oplus T^*#1}
\newcommand{\ttsc}[1]{T_\bC#1\oplus T_\bC^*#1}

\title[Local structure of generalized complex manifolds]{Local structure of
generalized complex manifolds}

\author[M.~Abouzaid and M.~Boyarchenko]{Mohammed
Abouzaid\address{Mohammed Abouzaid: Department of Mathematics,
University of Chicago, Chicago, IL 60637, e-mail:
mabouzai@math.uchicago.edu} \and Mitya Boyarchenko\address{Mitya
Boyarchenko: Department of Mathematics, University of Chicago,
Chicago, IL 60637, e-mail: mitya@math.uchicago.edu}}

\date{November 29, 2004}

\begin{document}

\begin{abstract}
We study generalized complex manifolds from the point of view of
symplectic and Poisson geometry. We start by showing that every
generalized complex manifold admits a canonical Poisson structure.
We use this fact, together with Weinstein's classical result on
the local normal form of Poisson manifolds, to prove a local
structure theorem for generalized complex manifolds which extends
the result Gualtieri has obtained in the ``regular'' case.
Finally, we begin a study of the local structure of a generalized
complex manifold in a neighborhood of a point where the associated
Poisson tensor vanishes. In particular, we show that in such a
neighborhood, a ``first-order approximation'' to the generalized
complex structure is encoded in the data of a constant $B$-field
and a complex Lie algebra.
\end{abstract}

\maketitle

\setcounter{tocdepth}{1}

\tableofcontents

\section{Introduction and main results}\label{s:intro}

The main objects of study in this paper are irregular generalized
complex (GC) structures on manifolds (the terminology is explained
below). In this section we state and discuss our main results. The
rest of the paper is devoted to their proofs.

\subsection{Background on GC geometry}\label{ss:background} We
begin by recalling the setup of generalized complex geometry. We
use \cite{marco} as the main source for most basic results and
definitions, a notable exception being the notion of a generalized
complex submanifold of a GC manifold, which is taken from
\cite{oren}.

\sbr

The notion of a GC manifold was introduced by N.~Hitchin (cf.
\cite{H1,H2,H3}) and developed by M.~Gualtieri in \cite{marco}. If
$M$ is a manifold (by which we mean a finite dimensional real
$C^\infty$ manifold), specifying a GC structure on $M$ amounts to
specifying either of the following two objects:
\begin{itemize}
\item an $\bR$-linear bundle automorphism $\cJ$ of $TM\oplus T^*M$
which preserves the standard symmetric bilinear pairing $\pair{
(X,\xi),(Y,\eta) } = \xi(Y)+\eta(X)$ and satisfies $\cJ^2=-1$, or
\item a complex vector subbundle $L\subset T_{\bC}M\oplus
T_{\bC}^*M$ such that $T_{\bC}M\oplus T_{\bC}^*M=L\oplus\bar{L}$
and $L$ is isotropic with respect to the $\bC$-bilinear extension
of $\pair{\cdot,\cdot}$ to $T_{\bC}M\oplus T_{\bC}^*M$,
\end{itemize}
which are required to satisfy a certain integrability condition
that is similar to the standard integrability condition for almost
complex structures on real manifolds. A bijection between the two
types of structure defined above is obtained by associating to an
automorphism $\cJ$ its $+i$-eigenbundle. In terms of $L$, the
integrability condition is that the sheaf of sections of $L$ is
closed under the {\em Courant bracket} \cite{cou}
\begin{equation}\label{e:courbra}
\cou{(X,\xi),(Y,\eta)}=\Bigl([X,Y],\,\cL_X\eta-\cL_Y\xi-\frac{1}{2}\cdot
d\bigl(\iota_X\eta-\iota_Y\xi\bigr)\Bigr).
\end{equation}
One can check that this condition is equivalent to the vanishing
of the {\em Courant-Nijenhuis tensor}
\begin{equation}\label{e:cntensor}
\cN_\cJ(A,B)=\cou{\cJ A,\cJ B} - \cJ\cou{\cJ A,B} - \cJ\cou{A,\cJ
B} - \cou{A,B},
\end{equation}
where $A,B$ are sections of $TM\oplus T^*M$.

\sbr

The two main examples of GC structures arise from complex and
symplectic manifolds. If $M$ is a real manifold equipped with an
integrable almost complex structure $J:TM\to TM$, it is easy to
check that the automorphism
\[
\cJ_J=\matr{J}{0}{0}{-J^*}
\]
defines a GC structure on $M$; such a GC structure is said to be
{\em complex}. Similarly, if $\om$ is a symplectic form on $M$, we
can view it as a skew-symmetric map $\om:TM\to T^*M$, and then the
automorphism
\[
\cJ_\om=\matr{0}{-\om^{-1}}{\om}{0}
\]
also defines a GC structure on $M$; such a GC structure is said to
be {\em symplectic}.

\sbr

A GC structure on a manifold $M$ induces a distribution
$E\subseteq T_\bC M$ which is smooth in the sense of \cite{Su}.
Namely, $E$ is the image of $L$ under the projection map
$\ttsc{M}\to\tc{M}$. Note that $E$ may not have constant rank. The
sheaf of sections of $E$ is closed under the Lie bracket (i.e.,
$E$ is involutive), as follows trivially from the definition of
the Courant bracket. Moreover, there is a (complex) $2$-form
$\eps$ on $E$ defined as follows: if $X,Y$ are sections of $E$,
choose a section $\xi$ of $\tcs{M}$ such that $(X,\xi)\in L$, and
set $\eps(X,Y)=\xi(Y)$. If $\eta$ is a section of $\tcs{M}$ such
that $(Y,\eta)\in L$, then $\xi(Y)=-\eta(X)$ because $L$ is
isotropic with respect to the pairing $\pair{\cdot,\cdot}$, which
implies that $\eps(X,Y)$ is independent of the choice of $\xi$;
thus $\eps$ is well defined. Furthermore, one can define the
tensor $d\eps\in\bigwedge^3(E^*)$ by the Cartan formula, which
makes sense since $E$ is involutive.

\begin{prop}[See \cite{marco}]\label{p:eps} The data $(E,\eps)$ determines the
GC structure $L$ uniquely. Moreover, $d\eps=0$.
\end{prop}

\sbr

A special type of operation defined for GC structures, which plays
an important role in our discussion, is the {\em transformation by
a $B$-field}. Specifically, if $B$ is a {\em real closed} $2$-form
on $M$, we define an orthogonal automorphism $\exp(B)$ of the
bundle $\tts{M}$ via
\[
\exp(B) = \matr{1}{0}{B}{1},
\]
where we view $B$ as a skew-symmetric map $TM\to T^*M$. If $\cJ$
defines a GC structure on $M$, and the associated pair $(E,\eps)$
is constructed as above, then $\cJ'=\exp(B)\cJ\exp(-B)$ is another
GC structure on $M$, which follows from the fact that $\exp(B)$
preserves the Courant bracket on $\tts{M}$, see \cite{marco}.
Moreover, in this case, the $+i$-eigenbundle of $\cJ'$ is given by
$L'=\exp(B)(L)$, and the associated pair $(E',\eps')$ is
determined by $E'=E$, $\eps'=\eps+B\big\vert_E$, where, by a
slight abuse of notation, we also denote by $B$ the $\bC$-bilinear
extension of $B$ to $\tc{M}$. In our paper, a $B$-field
transformation will always mean a transformation of the form
$\cJ\mapsto\exp(B)\cJ\exp(-B)$, where $B$ is a closed real
$2$-form. For a more detailed discussion and a more general notion
of $B$-fields, see \cite{marco,H3} and references therein.

\sbr

Another important construction is that of the canonical symplectic
foliation on a GC manifold. Namely, let us consider
$E\cap\bar{E}$; this is a distribution in $\tc{M}$ which is stable
under complex conjugation, and hence has the form
$\sS_\bC=\bC\otimes_\bR\sS$ for some distribution $\sS\subseteq
TM$. Gualtieri proves in \cite{marco} that $\sS$ is a smooth
distribution in the sense of \cite{Su}, and that the $2$-form
$\om$ on $\sS$ defined by $\om=\im\bigl(\eps\big\vert_\sS\bigr)$
is (pointwise) nondegenerate. Moreover, it is now clear that the
sheaf of sections of $\sS$ is closed under the Lie bracket, and
that $\om$ is a closed $2$-form on $\sS$, in the same sense as in
Proposition \ref{p:eps}. It follows from the results of \cite{Su}
that through every point of $M$ there is a maximal integral
manifold of $\sS$, which, by construction, inherits a natural
symplectic structure.

\sbr

For example, if the GC structure on $M$ is complex, then $\sS=0$,
while if the GC structure on $M$ is symplectic, then $\sS=TM$ and
the canonical symplectic form on $\sS$ coincides with the
symplectic form defining the GC structure on $M$.

\sbr

We now recall the notion of a generalized complex submanifold of a
GC manifold. Let $L$ be a GC structure on a manifold $M$, and let
$N\subset M$ be a (locally closed) submanifold. We define a (not
necessarily smooth) distribution $L_N$ on $N$ as follows. Set
\[
\tilde{L}_N=L\res{N}\cap \Bigl( \tc{N}\oplus
\bigl(\tcs{M}\res{N}\bigr) \Bigr) \quad\text{and}\quad
L_N=\pr\bigl(\tilde{L}_N\bigr),
\]
where $\pr:\tc{N}\oplus\bigl(\tcs{M}\res{N}\bigr)\to\ttsc{N}$
denotes the natural projection map, $\bigl(X,\xi\bigr)\mapsto
\bigl(X,\xi\big\vert_{\tc{N}}\bigr)$. It is proved in \cite{oren}
that $\dim_\bC L_{N,n}=\dim_\bR N$ for all $n\in N$. However,
$L_N$ may {\em not} be a subbundle of $\ttsc{N}$. We say that $N$
is a {\em generalized complex submanifold} of $M$ provided $L_N$
is smooth, and defines a GC structure on $N$. It can be shown (cf.
\cite{oren}) that a necessary and sufficient condition for this is
that $L_N$ is smooth and $L_N\cap\overline{L_N}=0$ (integrability
is then automatic).

\sbr

In conclusion, we would like to mention that there exists a way of
describing GC structures on manifolds in terms of spinors. In
fact, most of \cite{marco} is written in the language of spinors.
However, in our paper we have made a conscious effort to state and
prove all of our results in a spinor-free language. We hope that
this approach helps illuminate the simple geometric ideas that
underlie our main constructions.

\subsection{The canonical Poisson structure on a GC
manifold}\label{ss:poisson} From now on we fix a manifold $M$
equipped with a GC structure which, whenever convenient, we will
think of in terms of either the automorphism $\cJ$ or the
subbundle $L\subset\ttsc{M}$. The starting point for our work is
the observation that the canonical symplectic foliation
$(\sS,\om)$ defined in \S\ref{ss:background} is in fact the
symplectic foliation associated to a certain Poisson structure on
$M$. The existence of a canonical Poisson structure on a GC
manifold was also independently noticed by M.~Gualtieri
\cite{marco2}, and S.~Lyakhovich and M.~Zabzine \cite{lz}.

\sbr

Let us briefly explain why one could expect the existence of a
natural Poisson structure on general grounds. Recall the
definition of integrability as the vanishing of the
Courant-Nijenhuis tensor \eqref{e:cntensor}. The condition
$\cN_\cJ(A,B)=0$ can be naturally rewritten as a collection of
four equations corresponding to the possibilities of either $A$ or
$B$ being a section of $TM$ or a section of $T^*M$. Let us also
write $\cJ$ as a matrix
\begin{equation}\label{e:jmatrix}
\cJ=\matr{J}{\pi}{\sg}{K},
\end{equation}
where $J:TM\to TM$, $\pi:T^*M\to TM$, $\sg:TM\to T^*M$ and
$K:T^*M\to T^*M$ are bundle morphisms. The requirements that
$\cJ^2=-1$ and $\cJ$ is orthogonal with respect to
$\pair{\cdot,\cdot}$ force $K=-J^*$, $\pi=-\pi^*$, $\sg=-\sg^*$;
in particular, $\pi$ can be viewed as a bivector on $M$, i.e., a
section of $\bigwedge^2 TM$. Moreover, it is a straightforward
computation that in the case when $A=(0,\xi)$ and $B=(0,\eta)$,
where $\xi,\eta$ are sections of $T^*M$, the $TM$-component of
$\cN_\cJ(A,B)$ is the following expression:
\[
[\pi\xi,\pi\eta]-\pi\bigl( \cL_{\pi\xi}\eta-\frac{1}{2}
d(\iota_{\pi\xi}\eta) \bigr) + \pi\bigl(
\cL_{\pi\eta}\xi-\frac{1}{2} d(\iota_{\pi\eta}\xi) \bigr).
\]
Observe that this expression depends only on $\pi$ and not on the
other components of the matrix defining $\cJ$. However, one can
check that no other entry of the matrix can be separated from the
rest in this way. This suggests that $\pi$ must play a special
role in the theory. In fact, we prove
\begin{ithm}\label{t:poisson}
The bivector $\pi$ defines a Poisson structure on $M$. Moreover,
the canonical symplectic foliation associated to this Poisson
structure coincides with $(\sS,\om)$.
\end{ithm}
Given a real-valued $f\in C^\infty(M)$, let us write
\[
\cJ(0,df) = (X_f,\xi_f).
\]
By construction, $X_f=\pi(df)$ is the Hamiltonian vector field on
$M$ associated to $f$. On the other hand, $\xi_f$ is a certain
differential $1$-form on $M$.
\begin{prop}\label{p:xi} The map $f\mapsto\xi_f$ has the following
properties.
\begin{enumerate}
\item For all $f,g\in C^\infty(M)$, we have
\[
\xi_{f\cdot g} = f\cdot\xi_g + g\cdot\xi_f.
\]
\item If $\{\cdot,\cdot\}$ is the Poisson bracket on $C^\infty(M)$
defined by $\pi$, then
\[
\xi_{\{f,g\}}=\cL_{X_f}(\xi_g)-\iota_{X_g}(d\xi_f).
\]
\item If $(E,\eps)$ is associated to the GC structure $\cJ$ as in
\S\ref{ss:background}, then for all $f\in C^\infty(M)$, we have
\[
\cL_{X_f}(\eps) = \bigl(d\xi_f\bigr)\res{E}.
\]
\end{enumerate}
\end{prop}
The two results above are proved in Section \ref{s:poisson}. The
properties of the map $f\mapsto\xi_f$ turn out to be crucial in
our proof of the local normal form for GC manifolds. Moreover,
these result raise the question of whether one can give an
explicit description of GC manifolds as Poisson manifolds equipped
with additional structure. In other words, consider a GC structure
on a manifold $M$ defined by the matrix \eqref{e:jmatrix}. By
Theorem \ref{t:poisson}, the pair $(M,\pi)$ is a Poisson manifold.
Then the problem is to describe, in the language of Poisson
geometry, the extra data on $(M,\pi)$ that needed to recover all
of $\cJ$. Part (2) of Proposition \ref{p:xi} is a first step in
this direction.

\subsection{The local structure theorem for GC
manifolds}\label{ss:localstructure} We say that a GC structure on
a manifold $M$ is {\em regular} if the distribution $\sS$
(equivalently, $E$) has locally constant rank. The structure is
said to be {\em irregular} otherwise. The original motivation for
our work came from trying to extend the local structure theorem
proved in \cite{marco} for regular GC structures to the irregular
case. Gualtieri proved that if $m\in M$ is a regular point of a
given GC structure on $M$ (i.e., the structure is regular in an
open neighborhood of $M$), then there exists a neighborhood $U$ of
$m$ in $M$ such that the induced GC structure on $U$ is a
$B$-field transform of the product of a symplectic GC manifold and
a complex GC manifold. However, it seems to be difficult to adapt
the method of \cite{marco} to the irregular situation. In
particular, it relies strongly on the ``complex Frobenius
theorem'' \cite{cfrob}, no irregular analogue of which is known to
us. On the other hand, a powerful tool that is available to us in
view of Theorem \ref{t:poisson} is Weinstein's local structure
theorem for Poisson manifolds \cite{alan}. Our approach has the
advantage that it uses neither the real nor the complex version of
the Frobenius theorem; nor, indeed, any nontrivial result from the
theory of partial differential equations.

\sbr

Let us fix a GC manifold $M$ and a point $m_0\in M$. We define the
{\em rank}, $\rk_{m_0} M$, of $M$ at $m_0$ to be the rank of the
associated Poisson tensor $\pi$ at $m_0$. The central result of
our paper is the following
\begin{ithm}\label{t:local}
There exists an open neighborhood $U$ of $m_0$ in $M$, a real
closed $2$-form $B$ on $U$, a symplectic GC manifold $S$ and a GC
manifold $N$ with marked points $s_0\in S$, $n_0\in N$ such that
$\rk_{n_0}N=0$, and a diffeomorphism $S\times N\to U$ which takes
$(s_0,n_0)$ to $m_0$ and induces an isomorphism between the
product GC structure on $S\times N$ and the transform of the
induced GC structure on $U$ via the $2$-form $B$.
\end{ithm}

This theorem is proved in Section \ref{s:local}. Note that it is
different in nature from the recent results of Dufour and Wade
\cite{dw}. Due to the presence of $B$-fields, which have no
analogue for Dirac structures, our work gives more complete
information on the local structure of irregular GC manifolds than
loc.\,cit. does for irregular Dirac structures. Our method of
proof is also essentially different.

\begin{rem}
It is easy to recover the result of Gualtieri from Theorem
\ref{t:local}. Namely, if, with the notation of the theorem, the
GC structure on $M$ is regular in a neighborhood of $m_0$, then
the rank of $N$ must be zero in a neighborhood of $n_0$. It then
follows by linear algebra that the GC structure on $N$ must be
$B$-complex in a neighborhood of $n_0$, and the fact that this
structure can be written as the transform of a complex structure
by a {\em closed} real $2$-form follows from the local vanishing
of Dolbeault cohomology (cf. \cite{marco}).
\end{rem}

\subsection{Linear GC structures}\label{ss:linearGCS} The term ``linear GC
structure'' should not be confused with the notion of a {\em
constant} GC structure on a real vector space discussed in Section
\ref{s:linalg}. Rather, it is used in the same way as the term
``linear Poisson structure'' is used to describe the canonical
Poisson structure on the dual space of a real Lie algebra.

\sbr

Recall that if $(M,\pi)$ is a Poisson manifold, and $m\in M$ is a
point at which the Poisson tensor $\pi$ vanishes, then a
``first-order approximation'' to $\pi$ at $m$ defines a real Lie
algebra of dimension $\dim M$. Canonically, this Lie algebra can
be identified with the quotient $\fg=\fm_m/\fm_m^2$, where $\fm_m$
denotes the ideal in the algebra of all real-valued $C^\infty$
functions on $M$ consisting of the functions that vanish at $m$.
Since $\pi$ vanishes at $m$, it is easy to check that $\fm_m$ is
stable under the Poisson bracket, and $\fm_m^2$ is an ideal of
$\fm_m$ in the sense of Lie algebras, and hence we obtain an
induced Lie algebra bracket on $\fg$.

\sbr

Therefore one expects that, near a point on a GC manifold where
the associated Poisson tensor vanishes, the first-order
approximation to the GC structure can be encoded in a real finite
dimensional Lie algebra equipped with additional structure.
Indeed, we prove the following
\begin{ithm}\label{t:approx}
In a neighborhood of a point on a GC manifold where the associated
Poisson tensor vanishes, the first-order approximation to the GC
structure is encoded in a {\em complex} Lie algebra of complex
dimension $\bigl(\dim M\bigr)/2$, and a $B$-field which is
constant in appropriate local coordinates (and hence, a fortiori,
is closed).
\end{ithm}
The meaning of this statement is explained in Section
\ref{s:linear}.

\sbr

A natural problem that arises is to give a local classification of
GC manifolds near a point where the associated Poisson tensor
vanishes. Together with our Theorem \ref{t:local}, a solution of
this problem would yield a complete local classification of
generalized complex manifolds.

\subsection{Acknowledgements and credits} We would like to thank
Paul Seidel for his comments on the early versions of our paper,
and Marco Gualtieri for helpful discussions about $B$-fields.
After the first version of our paper was prepared, we learned from
Marco Gualtieri that the existence of a canonical Poisson
structure on generalized complex manifolds has been known to
physicists working in related areas, among them Lindstrom,
Lyakhovich, Minasian, Tomasiello, and Zabzine \cite{gcmsusy,lz},
and was explored by Crainic in a preprint \cite{cra}.
Nevertheless, as no proof seems to appear in the mathematical
literature, and since the result is central to our work, we
include a discussion of this Poisson structure and some of its
important properties.

\section{Linear algebra}\label{s:linalg}

\subsection{}\label{ss:linalgbasics} In this section we present the auxiliary results on
linear algebra that are used in the proofs of our main theorems.
We begin by recalling that the notion of a GC structure has an
analogue for vector spaces, which was studied in detail in
\cite{oren} and \cite{marco}. Specifically, a {\em constant}
generalized complex structure on a real vector space $V$ is
defined either as an $\bR$-linear automorphism $\cJ$ of $V\oplus
V^*$ which preserves the standard symmetric bilinear pairing
$\pair{\cdot,\cdot}$ and satisfies $\cJ^2=-1$, or as a complex
subspace $L\subset V_\bC\oplus V^*_\bC$ which is isotropic with
respect to the $\bC$-bilinear extension of $\pair{\cdot,\cdot}$
and satisfies $V_\bC\oplus V_\bC^*=L\oplus\Lb$. There is no
integrability condition in this case. It is easy to see that
constant GC structures on $V$ correspond bijectively to GC
structures on the underlying real manifold of $V$ that are
invariant under translations. Furthermore, it is obvious that if
$\cJ$ is a GC structure on a manifold $M$, then for every point
$m\in M$, the automorphism $\cJ_m$ of $T_m M\oplus T_m^* M$
induced by $\cJ$ defines a constant GC structure on $T_m M$. From
now on, by a {\em generalized complex vector space} we will mean a
real vector space equipped with a constant GC structure.

\sbr

All notions and constructions discussed in \S\ref{ss:background}
have obvious analogues for GC vector spaces. In particular, for a
real vector space $V$, we let $\rho:V\oplus V^*\to V$,
$\rho^*:V\oplus V^*\to V^*$ denote the natural projection maps.
Given a GC structure on $V$ defined by a subspace $L\subset
V_\bC\oplus V_\bC^*$, we let $E=\rho(L)\subseteq V_\bC$. There is
an induced $\bC$-bilinear $2$-form $\eps$ on $E$ defined in the
same way as in \S\ref{ss:background}, and the pair $(E,\eps)$
determines the GC structure on $V$ uniquely. Moreover, if
$S\subseteq V$ is the real subspace satisfying $\bC\otimes_\bR
S=E\cap\Eb$, then $\om=\im\bigl(\eps\big\vert_S\bigr)$ is a
symplectic form on $S$. Finally, the notion of a generalized
complex subspace of a GC vector space $V$ is defined in the
obvious way: if $W\subseteq V$ is a real subspace, set
\[
\tilde{L}_W=L\cap \bigl(W_\bC\oplus V_\bC^*\bigr)
\quad\text{and}\quad L_W=\pr(\tilde{L}_W),
\]
where $\pr:W_\bC\oplus V_\bC^*\to W_\bC\oplus W_\bC^*$ is the
projection map $(w,\la)\mapsto (w,\la\big\vert_{W_\bC})$. We say
that $W$ is a {\em generalized complex subspace} of $V$ if
$L_W\cap \overline{L_W}=(0)$; it is shown in \cite{oren} that in
this case $L_W$ is automatically a GC structure on $W$, called the
{\em induced generalized complex structure}.

\sbr

The notion of a $B$-field transform is also defined in the obvious
way. If $B\in\bigwedge^2 V^*$ is a skew-symmetric bilinear form on
$V$, then the map
\[
\exp(B)=\matr{1}{0}{B}{1}
\]
is a linear automorphism of $V\oplus V^*$ which preserves the
standard pairing $\pair{\cdot,\cdot}$, and hence acts on constant
GC structures on $V$ via
\[
L\mapsto \exp(B)\cdot L, \quad\text{or}\quad
\cJ\mapsto\exp(B)\cdot\cJ\cdot\exp(-B).
\]
It is easy to check that, in terms of the pairs $(E,\eps)$, the
transformation above is given by
\[
(E,\eps)\longmapsto \bigl( E, \eps+B_\bC\big\vert_E \bigr),
\]
where $B_\bC$ is the unique $\bC$-bilinear extension of $B$ to
$V_\bC$.

\sbr

In what follows, we will occasionally need to consider GC
structures on different vector spaces at the same time. Therefore,
whenever a confusion may arise, we will use the notation
$L_V\subset V_\bC\oplus V_\bC^*$, $\cJ_V\in\Aut_\bR(V\oplus V^*)$,
$S_V\subseteq V$, $E_V\subseteq V_\bC$, etc., to denote the
objects $L,\cJ,S,E$, etc., that are associated to a given GC
structure on a vector space $V$.

\subsection{} For future use, we make explicit the notions of an
isomorphism and a product of GC structures. Given two real vector
spaces, $P$ and $Q$, equipped with GC structures $L_P$ and $L_Q$,
an {\em isomorphism of GC vector spaces} between $P$ and $Q$ is an
$\bR$-linear isomorphism $\phi:P\to Q$ such that the induced map
\[
\bigl( \phi_\bC, (\phi_\bC^*)^{-1} \bigr) : P_\bC\oplus P_\bC^*
\rar{} Q_\bC \oplus Q_\bC^*
\]
carries $L_P$ onto $L_Q$. The {\em direct sum} of the GC vector
spaces $P$ and $Q$ is the vector space $P\oplus Q$ equipped with
the GC structure $L_P\oplus L_Q$ (called the {\em product GC
structure}), where we have made the natural identification
\[
(P\oplus Q)_\bC \oplus (P\oplus Q)_\bC^* \cong P_\bC\oplus P_\bC^*
\oplus Q_\bC \oplus Q_\bC^*.
\]
Finally, if $V$ is a GC vector space and $P,Q\subseteq V$ are two
subspaces, we say that $V$ is the direct sum of $P$ and $Q$ {\em
as GC vector spaces} provided $P,Q$ are GC subspaces of $V$, and
if we equip $P$, $Q$ with the induced GC structures and $P\oplus
Q$ with the product GC structure, then the map $P\oplus Q\to V$
given by $(p,q)\mapsto p+q$ is an isomorphism of GC vector spaces.

\sbr

The notions of an isomorphism and a product of GC structures have
obvious extensions to GC manifolds, see \cite{oren}.

\subsection{}The main results of generalized complex linear algebra
that we need are summarized in the following
\begin{thm}\label{t:gclinalg}
Let $V$ be any GC vector space, and let $(S,\om)$ be defined as
above.
\begin{enumerate}[{(}a{)}]
\item The notion of being a GC subspace is transitive; in fact,
the following stronger statement holds: if $W_1\subseteq V$ is a
GC subspace and $W_2\subseteq W_1$ is any real subspace, then
$W_2$ is a GC subspace of $V$ if and only if it is a GC subspace
$W_1$ with respect to the induced GC structure on $W_1$.
\footnote{In general, however, GC subspaces do not behave well
with respect to taking sums and intersections.} Moreover, if this
is the case, then the induced GC structure on $W_2$ is the same in
both cases.
\item A subspace $W\subseteq V$ is a GC subspace if and only if
$W\cap S$ is a symplectic subspace of $S$ (in the sense that
$\om\big\vert_{W\cap S}$ is nondegenerate) and $W_\bC=(W_\bC\cap
E)+(W_\bC\cap\Eb)$.
\item In particular, $S$ itself is a GC subspace of $V$; the
induced GC structure on $S$ is $B$-symplectic, and moreover, $S$
is the largest GC subspace of $V$ with this property. The
underlying symplectic structure on $S$ is given by $\om$.
\item The notion of being a GC subspace is invariant under
$B$-field transformations of the GC structure on $V$.
\item\label{i:compl} If $W\subseteq V$ is a real subspace such that $W+S=V$ (the
sum is not necessarily direct), then $W$ is a GC subspace of $V$
if and only if $W\cap S$ is a symplectic subspace of $S$. In
particular, any subspace of $V$ that is complementary to $S$ in
the sense of linear algebra is automatically a GC subspace of $V$.
\item\label{i:sum} Let $W\subseteq V$ be a real subspace such that $W+S=V$, and
let $S_0$ denote any real subspace of $S$ such that $S=S_0\oplus
(S\cap W)$, so that $V=S_0\oplus W$. Then the following two
conditions are equivalent:
\begin{enumerate}[(i)]
\item $S_0$ and $S\cap W$ are orthogonal with respect to $\om$;
\item $W$ and $S_0$ are GC subspaces of $V$, and there exists a
$B$-field $B\in\bigwedge^2 V^*$ which transforms the GC structure
on $V$ into the direct sum of the induced GC structures on $S_0$
and $W$.
\end{enumerate}
\item If the equivalent conditions of part (\ref{i:sum}) hold,
then the choice of $B$ is unique provided we insist that
$B\big\vert_{S_0}=0$ and $B\big\vert_W=0$.
\end{enumerate}
\end{thm}

\begin{rem}
As a byproduct of our discussion, we obtain an alternate proof of
the structure theorem for constant GC structures (see \cite{oren}
and \cite{marco}) which does not use spinors. Indeed, if
$S\subseteq V$ is as above and $W\subseteq V$ is any complementary
subspace to $S$, then parts (\ref{i:compl}) and (\ref{i:sum}) of
the theorem imply that $W$ is a GC subspace of $V$ and the GC
structure on $V$ is a $B$-field transform of the direct product GC
structure on $S\oplus W$. It is then easy to check that the
induced GC structure on $S$ (resp., $W$) is $B$-symplectic (resp.,
$B$-complex), see, e.g., \cite{oren}.
\end{rem}

\sbr

\begin{proof}[Proof of Theorem \ref{t:gclinalg}] \textbf{(a)} It is trivial
to check that the two definitions of $L_{W_2}$ we obtain by
viewing $W_2$ either as a subspace of $V$ or as a subspace of
$W_1$ coincide, whence the claim.

\sbr

\noindent
\textbf{(b)} We first show the necessity of the two conditions. It
follows from the results of \cite{oren} that a subspace of $S$ is
a GC subspace if and only if it is a symplectic subspace with
respect to the form $\om$. Now if $W$ is any GC subspace of $V$,
then $W\cap S=S_W$, whence $W\cap S$ is a GC subspace of $W$ by
the results of \cite{oren}. By part (a), it follows that $W\cap S$
is also a GC subspace of $V$, and hence a GC subspace of $S$.

\sbr

Suppose now that $W$ is a GC subspace of $V$, yet $(W_\bC\cap
E)+(W_\bC\cap\Eb)\subsetneq W_\bC$. Then there exists a nonzero
real subspace $U\subset W$ with
\[
U_\bC \oplus \bigl[ (W_\bC\cap E)+(W_\bC\cap\Eb) \bigr] = W_\bC.
\]
This implies that
\[
U_\bC\cap \bigl[E+(W_\bC\cap\Eb) \bigr] = (0)  \quad \text{and}
\quad U_\bC\cap \bigl[\Eb+(W_\bC\cap E) \bigr] = (0).
\]
Hence we can find $\ell,\ell'\in V_\bC^*$ with
$\ell\big\vert_{U_\bC}=\ell'\big\vert_{U_\bC}\not\equiv 0$ and
$\ell\big\vert_{E+(W_\bC\cap\Eb)}\equiv 0\equiv
\ell'\big\vert_{\Eb+(W_\bC\cap E)}$. This forces $\ell\in L\cap
V_\bC^*$, $\ell'\in \Lb\cap V_\bC^*$ and
$\ell\big\vert_{W_\bC}=\ell'\big\vert_{W_\bC}\not\equiv 0$, which
means that
\[
\bigl(\rho(\ell), \rho^*(\ell)\big\vert_{W_\bC} \bigr) = \bigl( 0,
\ell\big\vert_{W_\bC} \bigr) = \bigl( 0, \ell'\big\vert_{W_\bC}
\bigr) = \bigl(\rho(\ell'), \rho^*(\ell')\big\vert_{W_\bC} \bigr)
\neq 0,
\]
contradicting the assumption that $W$ is a GC subspace of $V$.

\sbr

Conversely, suppose that $W\subseteq V$ is a subspace such that
$W_\bC=(W_\bC\cap E)+(W_\bC\cap\Eb)$ and $W\cap S$ is a GC
subspace (equivalently, a symplectic subspace) of $S$. We will
prove that $W$ is a GC subspace of $V$. Assume that $\ell\in L$,
$\ell'\in\Lb$ and $\rho(\ell)=\rho(\ell')\in W_\bC$,
$\rho^*(\ell)\big\vert_{W_\bC}=\rho^*(\ell')\big\vert_{W_\bC}$.
Then, in particular, $\rho(\ell)=\rho(\ell')\in (W\cap S)_\bC$ and
$\rho^*(\ell)\big\vert_{(W\cap
S)_\bC}=\rho^*(\ell')\big\vert_{(W\cap S)_\bC}$, so we deduce from
the second assumption that $\rho(\ell)=\rho(\ell')=0$ and
$\rho^*(\ell)\big\vert_{(W\cap
S)_\bC}=\rho^*(\ell')\big\vert_{(W\cap S)_\bC}=0$. It remains to
check that
$\rho^*(\ell)\big\vert_{W_\bC}=\rho^*(\ell')\big\vert_{W_\bC}=0$.
But
\[
\rho^*(\ell)=\ell\in L\cap V_\bC^* = \Ann_{V_\bC^*}(E) \quad
\text{and} \quad \rho^*(\ell')=\ell'\in \Lb\cap V_\bC^* =
\Ann_{V_\bC^*}(\Eb),
\]
whence $\ell\big\vert_{W_\bC\cap E}=0=\ell'\big\vert_{W_\bC\cap
\Eb}$, and also, since
$\ell\big\vert_{W_\bC}=\ell'\big\vert_{W_\bC}$, we find from our
first assumption that
$\ell\big\vert_{W_\bC}=\ell'\big\vert_{W_\bC}=0$, completing the
proof.

\sbr

\noindent
\textbf{(c)} This is easy. We omit the proof since the
straightforward argument is presented in \cite{oren}.

\sbr

\noindent
\textbf{(d)} It follows from the remarks of
\S\ref{ss:linalgbasics} that a $B$-field transform changes neither
$E$, nor $S$, nor $\om=\im\bigl(\eps\big\vert_S\bigr)$. Hence the
claim follows from the characterization of GC subspaces given in
part (b).

\sbr

\noindent
\textbf{(e)} We will show that if $W\subseteq V$ is a subspace
such that $V=W+S$, then we automatically have $W_\bC=(W_\bC\cap
E)+(W_\bC\cap\Eb)$. The claim then follows from part (b). Let
$w\in W_\bC$, and write $w=e_1+\overline{e_2}$, with $e_j\in E$
for $j=1,2$. Further, we can write $e_j=w_j+s_j$, where $w_j\in
W_\bC$ and $s_j\in S_\bC$. A fortiori, $s_j\in E$, so $w_j\in
E\cap W_\bC$. Hence
\[
w=w_1+\barr{w_2}+(s_1+\barr{s_2}),
\]
where $w_1,w_2\in W_\bC\cap E$. This forces $s_1+\barr{s_2}\in
W_\bC$, and since we also have $s_j\in S_\bC$, it follows that
$s_1+\barr{s_2}\in W_\bC\cap S_\bC\subseteq W_\bC\cap E$. Finally,
we conclude that
\[
w=(w_1+s_1+\barr{s_2})+\barr{w_2},
\]
where $w_1+s_1+\barr{s_2}\in W_\bC\cap E$ and $\barr{w_2}\in
W_\bC\cap\Eb$, as desired.

\sbr

\noindent
\textbf{(f), (g)} First, it is clear that (ii) implies (i), since
$B$-field transforms cannot change the imaginary part of $\eps$.
Conversely, assume that $S_0$ and $S\cap W$ are orthogonal with
respect to $\om$. We will show that there exists exactly one
$B$-field $B\in\bigwedge^2 V^*$ such that
$B\big\vert_{S_0}=B\big\vert_W=0$ and $B$ transforms the given GC
structure on $V$ into the direct sum of the induced GC structures
on $S_0$ and $W$.

\sbr

Observe that $E_V=E_{S_0}\oplus E_W$. Indeed, it is clear that
$E_{S_0}\oplus E_W\subseteq E_V$. Conversely, let $e\in E_V$ and
write $e=e_1+e_2$, where $e_1\in (S_0)_\bC$ and $e_2\in W_\bC$.
Then, a fortiori, $e_1\in E_V$, so we also have $e_2\in E_V$,
whence $e_1\in E_V\cap (S_0)_\bC=E_{S_0}$ and $e_2\in E_V\cap
W_\bC=E_W$, proving the claim.

\sbr

Note now that if the original GC structure on $V$ is determined by
$(E_V,\eps)$, then the product GC structure on $S_0\oplus W$ is
determined by
\[
\Bigl( E_{S_0}\oplus E_W, \eps\big\vert_{E_{S_0}} +
\eps\big\vert_{E_W} \Bigr).
\]
To complete the proof, we must therefore show that there exists
exactly one $B\in\bigwedge^2 V^*$ such that
$B\big\vert_{S_0}=B\big\vert_W=0$ and the pairing between
$E_{S_0}$ and $E_W$ induced by (the complexification of) $B$ is
the same as the one induced by $\eps$.

\sbr

Suppose that such a $B$ exists. Let $s\in S_0$, $w\in W$. Since
$w$ is real, we can write $w=e+\barr{e}$, where $e\in E\cap
W_\bC$. Then we must have
\[
B(s,w)=B(s,e)+\barr{B(s,e)}=2\cdot\re\eps(s,e),
\]
which proves that $B$ is unique if it exists. Conversely, let us
define $B$ on $S_0\times W$ by this formula, and define $B$ to be
zero on $S_0$ and on $W$. We claim that $B$ is well defined.
Indeed, consider a different representation $w=e'+\barr{e'}$,
where $e'\in E\cap W_\bC$. Then
\[
e-e' = \barr{e'-e} \in (W\cap S)_\bC,
\]
which implies that $e-e'=i\cdot t$ for some $t\in W\cap S$, where
$i=\sqrt{-1}$. Hence
\[
\re\eps(s,e-e') = -\im\eps(s,t) = -\om(s,t) = 0 \quad \text{by
assumption},
\]
which implies that $B$ is well defined.

\sbr

Finally, to show that $B$ satisfies the required condition, it is
enough to check (by linearity) that if $s\in S_0$ and $e\in
E_W=W_\bC\cap E$, then $B(s,e)=\eps(s,e)$. We have
\[
e= \frac{e+\barr{e}}{2} + i\cdot \frac{e-\barr{e}}{2i}
\quad\text{and}\quad \frac{e+\barr{e}}{2},
\frac{e-\barr{e}}{2i}\in W.
\]
By construction,
\[
B\left(s,\frac{e+\barr{e}}{2}\right) = 2\re\eps
\left(s,\frac{e}{2}\right) = \re\eps(s,e),
\]
and similarly $B\bigl(s,(e-\barr{e})/(2i)\bigr)=\im\eps(s,e)$,
which completes the proof.
\end{proof}

\begin{rem}\label{r:smoothness}
The following comment will be used in our proof of the local
structure theorem for GC manifolds. Consider a variation of
generalized complex linear algebra where the vector space $V$ is
replaced by a smooth real vector bundle $\cV$ over a base manifold
$\fB$, and a GC structure on $\cV$ is a subbundle
$\cL\subseteq\cV_{\bC}\oplus\cV^*_{\bC}$ such that for every point
$b\in\fB$, the subspace
$\cL_b\subseteq\cV_{b,\bC}\oplus\cV^*_{b,\bC}$ defines a constant
GC structure on the real vector space $\cV_b$. Then we have the
subdistributions $\cE\subseteq\cV_\bC$ and $\cS\subseteq\cV$ which
are the global analogues of $E$ and $S$, respectively, which may
have nonconstant rank, but are nevertheless smooth in the sense of
\cite{Su}, by the argument given in \cite{marco}. We claim that,
in fact, the proofs of parts (e), (f) and (g) of Theorem
\ref{t:gclinalg} go through in this setup with appropriate
modifications that ensure smooth dependence on the point
$b\in\fB$.

\sbr

First, consider the analogue of part (e), where $W$ is replaced by
a smooth subbundle $\cW\subseteq\cV$ such that $\cV=\cS+\cW$
pointwise. We assume also that there exists a smooth subbundle
$\cS_0\subseteq\cV$ which is contained in $\cS$ and satisfies
$\cV=\cS_0\oplus\cW$. Then we claim that
$\cW_\bC=(\cW_\bC\cap\cE)+(\cW_\bC\cap\barr{\cE})$ in the sense
that every smooth section $w$ of $\cW_\bC$ can be written as
$w=w'+w''$, where $w'$ and $w''$ are smooth sections of $\cW_\bC$
that lie in $\cE$ and $\barr{\cE}$, respectively. Indeed, let us
go through the proof of part (e) given above. By assumption,
$\cV_\bC\oplus\cV_\bC^*=\cL\oplus\barr{\cL}$, so we can write
$(w,0)=(e_1,f_1)+(e_2,f_2)$ for smooth sections $(e_1,f_1)$ and
$(e_2,f_2)$ of $\cL$ and $\barr{\cL}$, respectively. A fortiori,
$e_1$ and $e_2$ are smooth sections of $\cE$ and $\barr{\cE}$ that
satisfy $w=e_1+e_2$. Further, we can write uniquely $e_j=w_j+s_j$,
where $w_j$ are smooth sections of $\cW_\bC$ and $s_j$ are smooth
sections of $\cS_{0,\bC}$. The rest of the proof of part (e) goes
through without changes.

\sbr

Next we consider the analogue of parts (f) and (g). We assume that
$\cS_0$, $\cW$ are as in the previous paragraph, and that, in
addition, $\cS_0$ and $\cW$ are orthogonal with respect to the
canonical symplectic form on $\cS$. If $B$ is the $2$-form on
$\cV$ constructed pointwise as in the proof of parts (f) and (g)
given above, we claim that $B$ is in fact smooth. Clearly, it
suffices to show that if $s$ and $w$ are smooth sections of
$\cS_0$ and $\cW$, respectively, then $B(s,w)$ is a smooth
function on $\fB$. By the previous paragraph, we can write
$w=e+\barr{e'}$, where $e,e'$ are smooth sections of $\cW_\bC$
that lie in $\cE$. Since $w$ itself is real, we have
\[
w=\frac{1}{2}\cdot \bigl[ (e+\barr{e'})+\barr{(e+\barr{e'})}
\bigr] = \frac{1}{2}\cdot\bigl[ (e+e')+\barr{e+e'} \bigr],
\]
which implies that we may assume that $e=e'$ without sacrificing
smoothness. By the proof above, we then have
$B(s,w)=2\cdot\re\eps(s,e)$, which is a smooth function since
$\eps$ is smooth.
\end{rem}

\subsection{} As one can see from part (\ref{i:sum}) of Theorem
\ref{t:gclinalg}, if $W\subseteq V$ is a GC subspace such that
$S+W=V$, it is important to know the orthogonal complement $(W\cap
S)^\perp$ of $W\cap S$ in $S$ with respect to $\om$. The following
useful result provides a construction of this orthogonal
complement that does not involve $S$ or $\om$ explicitly.

\begin{thm}\label{t:orthcompl}
Let us identify $\Lb$ with $L^*$ using the standard symmetric
bilinear pairing $\pair{\cdot,\cdot}$, and define $A_W$ to be the
annihilator inside $\Lb$ of the subspace $\tilde{L}_W=L\cap
(W_\bC\oplus V_\bC^*)$. If $C_W$ is the projection of $A_W$ on
$V_\bC$, then $(W\cap S)^\perp_\bC=C_W$.
\end{thm}

The proof of this theorem provides an illustration of the
techniques developed in this section. It consists of four steps.

\sbr

\noindent
\textbf{Step 1.} The statement of the theorem holds when the GC
structure on $V$ is symplectic, given by a symplectic form
$\om\in\bigwedge^2 V^*$.
\begin{proof}
In this case $V=S$, the canonical symplectic form on $S$ is also
given by $\om$, and the condition $W+S=V$ is vacuous. We have
\[
L=L_V=\bigl\{ (v,-i\cdot\om(v)) \big\vert v\in V_\bC \bigr\},
\]
whence
\[
\tilde{L}_W=\bigl\{ (w,-i\cdot\om(w)) \big\vert w\in W_\bC
\bigr\}.
\]
Also,
\[
\Lb = \bigl\{ (u,i\cdot\om(u)) \big\vert u\in V_\bC \bigr\},
\]
whence
\begin{eqnarray*}
A_W &=& \Bigl\{ \bigl( u,i\om(u) \bigr) \Big\vert u\in V_\bC,\
i\om(u,w)-i\om(w,u)=0 \text{ for all } w\in W_\bC \Bigr\} \\
&=& \Bigl\{ \bigl( u,i\om(u) \bigr) \Big\vert u\in V_\bC,\
\om(u,w)=0 \text{ for all } w\in W_\bC \Bigr\},
\end{eqnarray*}
and so
\[
C_W = \Bigl(W^{\perp{\om}}\Bigr)_\bC,
\]
as required.
\end{proof}

\sbr

\noindent
\textbf{Step 2.} The statement of the theorem holds under the
following assumption: there exists a GC subspace $U\subseteq V$
such that
\begin{itemize}
\item $W=(W\cap S)\oplus U$;
\item the induced GC structure on $S$ (resp., $U$) is symplectic
(resp., complex);
\item $V$ is the direct sum of $S$ and $U$ {\em as GC vector
spaces}.
\end{itemize}
\begin{proof}
We make the obvious identifications
\[
V^*\cong S^*\oplus U^*,\quad V_\bC^*\cong S_\bC^*\oplus U_\bC^*,
\quad V_\bC\oplus V^*_\bC\cong S_\bC\oplus S_\bC^* \oplus U_\bC
\oplus U_\bC^*.
\]
Under these identifications, we have, by assumption,
\[
L=L_V=L_S\oplus L_W \quad\text{and}\quad
\tilde{L}_W=\tilde{L}_{(W\cap S)}\oplus L_U.
\]
Thus it is clear that
\[
A_W=A_{W\cap S}\oplus (0)\subseteq \Lb_S\oplus\Lb_U=\Lb_V,
\]
whence the claim follows from Step 1.
\end{proof}

\sbr

\noindent
\textbf{Step 3.} The statement of the theorem is invariant under
$B$-field transformations.
\begin{proof}
Let $B\in\bigwedge^2 V^*$, then (with self-explanatory notation)
we have
\[
L^{new} = \exp(B)(L^{old}), \quad \Lb^{new} = \exp(B)(\Lb^{old}),
\]
\[
\tilde{L}_W^{new} = \exp(B)(\tilde{L}_W^{old}) \quad
(\text{because } \exp(B) \text{ preserves the subspace } W\oplus
V^*\subseteq V\oplus V^*),
\]
\[
A_W^{new} = \exp(B)(A_W^{old}) \quad (\text{because } \exp(B)
\text{ is orthogonal with respect to the standard pairing }
\pair{\cdot,\cdot}),
\]
and finally
\[
C_W^{new} = C_W^{old} \quad (\text{because } \exp(B) \text{
commutes with the projection onto } V_\bC).
\]
On the other hand, we know that a $B$-field transformation changes
neither $S$ nor the canonical symplectic form on $S$, proving the
claim.
\end{proof}

\sbr

\noindent
\textbf{Step 4.} We now complete the proof of the theorem as
follows. Since $W+S=V$, there exists a subspace $U\subseteq W$
such that $V=S\oplus U$. It follows from Theorem \ref{t:gclinalg}
that $U$ is a GC subspace of $V$, and the GC structure on $V$ is a
$B$-field transform of the direct sum of the induced GC structures
on $U$ and $S$. Moreover, it follows from the results of
\cite{oren} that the induced GC structure on $S$ (resp., $U$) is
$B$-symplectic (resp., $B$-complex). Hence there exists a
$B\in\bigwedge^2 V^*$ that transforms the GC structure on $V$ into
the direct sum of the underlying symplectic structure on $S$ and
the underlying complex structure on $U$. Now Steps 2 and 3
complete the argument.

\section{Generalized complex manifolds are
Poisson}\label{s:poisson}

In this section we prove Theorem \ref{t:poisson} and Proposition
\ref{p:xi}. However, it is convenient for us to first restate the
definition of the Poisson bracket on a generalized complex
manifold in a different way.

\sbr

Let $M$ be a manifold and $C^\infty(M)$ the algebra of all
real-valued smooth functions on $M$. Consider a GC structure on
$M$, where $\cJ$ is the corresponding automorphism of $\tts{M}$
and $L$ is the $+i$-eigenbundle of $\cJ$. If $f\in C^\infty(M)$,
then $(0,df)\in\tts{M}$ can be written as
$(0,df)=(X,\xi)+(\overline{X'},\overline{\xi'})$ for some sections
$(X,\xi)$, $(X',\xi')$ of $L$. Because of the uniqueness of this
decomposition, and since $df=\overline{df}$, we have
$(X',\xi')=(X,\xi)$. We set
\[
X_f = 2iX = -2\im(X),
\]
and call it the {\em Hamiltonian vector field} associated to the
function $f$. Further, we put $\xi_f=-2\im(\xi)$. Since $L$
(resp., $\Lb$) is the $+i$-eigenbundle (resp., $-i$-eigenbundle)
of $\cJ$, it is clear that
\[
\cJ(0,df) = (X_f,\xi_f).
\]
Recall that $(\sS,\om)$ denotes the canonical symplectic foliation
of $M$, as defined in \S\ref{ss:background}.
\begin{lem}\label{l:p1}
The vector field $X_f$ lies in $\sS$. Moreover, for every section
$Y$ of $\sS$, we have
\[
\om(X_f,Y)=Y(f).
\]
\end{lem}
\begin{proof}
We use the notation of the previous paragraph. Since
$X+\overline{X'}=X+\Xb=0$, we see that $X_f$ is a real vector
field. Moreover, $X_f=2iX$ lies in $E$ by construction, which
forces $X_f$ to lie in $\sS$. Now let $Y$ be a section of $\sS$.
By definition, we have
\[
\om(X_f,Y)=\im\eps(X_f,Y)=2\cdot\re\eps(X,Y)=2\cdot\re\xi(Y)=(df)(Y)=Y(f).
\]
\end{proof}
\begin{lem}\label{l:p2}
With the notation above, the flow of $X_f$ preserves the
subbundles $\sS\subseteq TM$, $E\subseteq\tc{M}$, and also
preserves the symplectic form $\om$ on $\sS$.
\end{lem}
\begin{proof}
The fact that the flow of $X_f$ preserves $\sS$ and $E$ follows
from the fact that $\sS\subset E$, that $X_f$ lies in $\sS$, and
that both $\sS$ and $E$ are integrable. Then the fact that $\om$
is preserved is a standard computation:
\[
\cL_{X_f}\om = d(\iota_{X_f}\om)+\iota_{X_f}(d\om)=d(df)+0=0.
\]
\end{proof}
\begin{lem}\label{l:p3}
Let $f,g\in C^\infty(M)$, and define $\{f,g\}=X_f(g)$. Then
\[
X_{\{f,g\}}=[X_f,X_g].
\]
\end{lem}
\begin{proof}
With the same notation as above, write
\[
(0,df)=(X,\xi)+(\Xb,\overline{\xi}), \quad
(0,dg)=(Y,\eta)+(\Yb,\overline{\eta}),
\]
so that
\[
X_f=2iX,\quad X_g=2iY.
\]
Recall from \cite{marco} that the restriction of the Courant
bracket to sections of $L$ can be written as follows:
\[
\cou{(X,\xi),(Y,\eta)}=\bigl( [X,Y],\cL_X\eta-\iota_Y(d\xi)
\bigr).
\]
Therefore the integrability condition for $L$ implies that
\begin{equation}\label{e:poiss}
\ell=\Bigl( \frac{1}{2i}[X_f,X_g], \cL_{X_f}\eta -
\iota_{X_g}(d\xi) \Bigr) = 2i\cdot \bigl( [X,Y],
\cL_X\eta-\iota_Y(d\xi) \bigr)
\end{equation}
is a section of $L$. Now the first component of
$\ell+\overline{\ell}$ is zero, and the second one is
\[
\cL_{X_f}(\eta+\overline{\eta})-\iota_{X_g}(d\xi+d\overline{\xi})
=
\cL_{X_f}(dg)-\iota_{X_g}(ddf)=\iota_{X_f}(ddg)+d(\iota_{X_f}(dg))=d
X_f(g)=d\{f,g\}.
\]
Consequently,
\[
X_{\{f,g\}}=2i\cdot \frac{1}{2i} [X_f,X_g]=[X_f,X_g].
\]
\end{proof}
\begin{thm}\label{t:poiss}
The bracket $\{\cdot,\cdot\}$ defined in Lemma \ref{l:p3} is a
Poisson bracket on $M$. Moreover, the Hamiltonian vector fields
$X_f$ locally span the distribution $\sS\subseteq TM$.
Consequently, $(\sS,\om)$ coincides with the canonical symplectic
foliation associated to the Poisson structure $\{\cdot,\cdot\}$.
\end{thm}
\begin{proof}
By Lemma \ref{l:p1}, we have $\{f,g\}=\om(X_g,X_f)$, which shows
that $\{\cdot,\cdot\}$ is skew-symmetric. The Leibniz rule is
straightforward:
\[
\{f,gh\}=X_f(gh)=g X_f(h)+ h X_f(g) = g\{f,h\}+h\{f,g\}.
\]
To check the Jacobi identity, we compute, using Lemma \ref{l:p3}:
\[
\{\{f,g\},h\}=X_{\{f,g\}}(h)=[X_f,X_g](h)=X_f X_g(h) - X_g X_f(h)
= \{f,\{g,h\}\}-\{g,\{f,h\}\}.
\]
This proves the first statement of the theorem. The last statement
follows immediately from the first two and the identity
$\{f,g\}=\om(X_g,X_f)$ proved in Lemma \ref{l:p1}. Thus it remains
to prove the second statement of the theorem. The question is a
pointwise one; thus let us choose a point $m\in M$, and let
$\{f_j\}$ be any collection of functions in $C^\infty(M)$ such
that their differentials $df_j$ span the cotangent space $T_m^*M$.
If $Y$ is any local section of $\sS$ near $m$ such that
$\om(Y,X_{f_j})=0$ at $m$ for all $j$, then Lemma \ref{l:p1}
implies that $\pair{Y,df_j}=0$ at $m$ for all $j$, whence $Y$
vanishes at $m$. This shows that the vector fields $X_{f_j}$ span
$\sS$ at $m$.
\end{proof}
Note that Theorem \ref{t:poisson} in the introduction is identical
to the result we have just proved. We can now give a
\begin{proof}[Proof of Proposition \ref{p:xi}]
Part (1) follows from the fact that $\xi_f$ is defined by the
equation $\cJ(0,df)=(X_f,\xi_f)$, that $(0,d(fg))=f\cdot
(0,dg)+g\cdot (0,df)$, and that $\cJ$ is linear over
$C^\infty(M)$. We prove (2) using the notation of the proof of
Lemma \ref{l:p3}. If $\ell$ is defined by \eqref{e:poiss}, then,
since $\ell$ is a section of $L$ and
$\ell+\overline{\ell}=(0,d\{f,g\})$, we have, by definition,
\[
\xi_{\{f,g\}}=-2\im\bigl( \cL_{X_f}\eta-\iota_{X_g}(d\xi) \bigr) =
\cL_{X_f}(-2\im\eta)-\iota_{X_g}\bigl(d(-2\im\xi)\bigr)=\cL_{X_f}(\xi_g)-\iota_{X_g}(d\xi_f).
\]
Finally, to prove (3), let us use the notation of the beginning of
this section and write $(0,df)=(X,\xi)+(\Xb,\overline{\xi})$,
where $(X,\xi)$ is a section of $L$. Then $X_f=2iX$, so
$(X_f,2i\xi)$ is also a section of $L$, which implies, by the
definition of $\eps$, that
\[
\iota_{X_f}(\eps)=2i\xi\res{E}=2i\cdot\bigl(\re\xi+i\im\xi)\res{E}=\bigl(i\cdot
df+\xi_f\bigr)\res{E}.
\]
Now we compute, using Cartan's formula and the fact that $d\eps=0$
(cf. Proposition \ref{p:eps}):
\[
\cL_{X_f}(\eps) = d\bigl(\iota_{X_f}(\eps)\bigr) +
\iota_{X_f}(d\eps) = d\bigl(\iota_{X_f}(\eps)\bigr) =
d\bigl(i\cdot df+\xi_f\bigr)\res{E} = \bigl(d\xi_f\bigr)\res{E}.
\]
This completes the proof.
\end{proof}

\section{Local normal form}\label{s:local}

\subsection{Strategy of the proof}\label{ss:strategy} We begin by outlining the
strategy of our proof of Theorem \ref{t:local}. Our argument is an
extension of the inductive argument of \cite{alan}. If
$\rk_{m_0}M=0$, then there is nothing to prove. Otherwise,
following loc.\,cit., we can split $M$, locally near $m_0$, as a
product $M=S\times N$ in the sense of Poisson manifolds,
$M=S\times N$, where $S$ is an open neighborhood of $0$ in $\bR^2$
with the induced standard symplectic form $\om_0$, and $m_0\in M$
corresponds to $(0,n_0)\in S\times N$. By abuse of notation, we
identify $N$ with the submanifold $\{0\}\times N$ of $M$. It is
clear that each ``horizontal leaf'' $S\times\{n\}$ is a GC
submanifold of $M$.
\begin{lem}\label{l:subman}
The ``transverse slice'' $N$ is a GC submanifold of $M$.
\end{lem}
The proof of this lemma is given at the end of the section. We
equip $N$ with the induced GC structure. It is clear that
$\rk_{n_0}N=\rk_{m_0}M-2$. Hence, by induction, it suffices to
show that in a neighborhood of $m_0$, the GC structure on $M$ is a
$B$-field transform of the product of the symplectic structure on
$S$ and the induced GC structure on $N$. The proof of this fact
consists of three steps, each involving a transformation by a
closed $2$-form and possibly replacing $M$ by a smaller open
neighborhood of $m_0$. To save space, we will still use $M$ to
denote any of these sufficiently small neighborhoods. The steps
are listed below:
\begin{enumerate}[(1)]
\item After a transformation by a closed $2$-form $B''$ on $M$, the
induced GC structure on each horizontal leaf $S\times\{n\}$ is the
symplectic GC structure defined by $\om_0$ via the obvious
identification $S\cong S\times\{n\}$.
\item After a transformation by a closed $2$-form $B'$ on $M$ that
restricts to zero on the horizontal leaves $S\times\{n\}$ and on
the transverse slice $\{0\}\times N$, we have that for each $n\in
N$, the induced constant GC structure on $T_{(0,n)}M$ is the
direct sum of the induced constant GC structures on
$T_{(0,n)}\bigl(S\times\{n\}\bigr)$ and on $T_n N$.
\item After a transformation by a closed $2$-form $B$ on $M$ that
vanishes along $N$, the GC structure on all of $M$ is the product
of the symplectic GC structure on $S$ and the induced GC structure
on $N$.
\end{enumerate}

\subsection{Step 1}\label{ss:step1} We begin by introducing
notation that will be used in the rest of the section. Let $(p,q)$
denote the standard coordinates on $S$, so that $\om_0=dp\wedge
dq$; we will also view them as part of a coordinate system
$(p,q,r_1,\dotsc,r_d)$ on $M$, where $r_1,\dotsc,r_d$ are local
coordinates on $N$ centered at $n_0$. Note that for any such
coordinate system on $M$, we have
\begin{equation}\label{e:xpxq}
X_p = -\frac{\dd}{\dd q} \quad\text{and}\quad X_q=\frac{\dd}{\dd
p},
\end{equation}
where $X_p$ and $X_q$ denote the Hamiltonian
vector fields on $M$ associated to the functions $p$ and $q$, as
in \S\ref{s:poisson}.

\sbr

Without loss of generality, we may assume that $S$ is the open
square on $\bR^2$ defined by the inequalities $-1<p<1$, $-1<q<1$.
A point $(s,n)\in S\times N=M$ will from now on be written as
$(a,b,n)$, where $a=p(s)$, $b=q(s)\in (-1,1)$.  We will denote by
$\phi_s:M\to M$ and $\psi_t:M\to M$ the flows of the vector field
$X_p$ and $X_q$, respectively. Of course, these flows are not
everywhere defined. Explicitly, we have, from \eqref{e:xpxq},
\begin{equation}\label{e:flows}
\phi_s(a,b,n)=(a,b-s,n) \quad\text{and}\quad
\psi_t(a,b,n)=(a+t,b,n).
\end{equation}
It is clear that the flows $\phi_s$ and $\psi_t$ commute with each
other.

\sbr

Furthermore, we define $\sS_0$ (resp., $\sN$) to be the
distribution on $M$ which is tangent to the horizontal leaves
$S\times\{n\}$ (resp., to the transverse slices $\{s\}\times N$);
note that $\sS_0$ is spanned by the vector fields $X_p$, $X_q$.

\sbr

We now prove statement (1) of \S\ref{ss:strategy}. Since
$M=S\times N$ as Poisson manifolds, it follows that for each $n\in
N$, the induced GC structure on $S\times\{n\}$ is $B$-symplectic,
with the underlying symplectic structure being given by $\om_0$. A
general fact, proved in \cite{marco}, is that on a $B$-symplectic
GC manifold, both the underlying symplectic structure and the
$B$-field are uniquely determined, and, moreover, depend smoothly
on the original GC structure. In our situation, we obtain a family
$\{B_n''\}$ of closed $2$-forms on the leaves $S\times\{n\}$,
depending smoothly on $n$, such that for every $n\in N$, the
$B$-field $B_n''$ transforms the induced GC structure on
$S\times\{n\}$ into the symplectic structure on $S\times\{n\}$
defined by $\om_0$.

\sbr

The usual proof of the Poincar\'e lemma shows that, after possibly
shrinking $S$ and $N$, we can find a smooth family
$\{\sg_n\}_{n\in N}$ of $1$-forms on the leaves $S\times\{n\}$
such that $d\sg_n=B_n''$ for each $n\in N$. Now let $\sg$ be an
arbitrary $1$-form on $M$ such that $\sg\res{S\times\{n\}}=\sg_n$
for each $n\in N$; such a $\sg$ exists simply because
$TM=\sS_0\oplus\sN$ as vector bundles. By construction, the
$2$-form $B''=d\sg$ satisfies the requirement of statement (1) of
\S\ref{ss:strategy}.

\subsection{Step 2}\label{ss:step2} It follows now from parts
(f) and (g) of Theorem \ref{t:gclinalg}, together with Remark
\ref{r:smoothness}, that for every point $n\in N$, there exists a
unique $2$-form $B_n'\in\Wedge^2 T^*_{(0,0,n)}M$ with the
following properties:
\begin{itemize}
\item $B_n'\res{T_n N}=0$;
\item $B_n'\res{T_{(0,0,n)}\bigl(S\times\{n\}\bigr)}=0$;
\item $B_n'$ transforms the constant GC structure on
$T_{(0,0,n)}M$ into the direct sum of the induced GC structures on
$T_n N$ and $T_{(0,0,n)}\bigl(S\times\{n\}\bigr)$;
\end{itemize}
and, moreover, $B_n'$ depends smoothly on $n$. We must show that
there exists a closed $2$-form $B'$ on $M$ such that for each
$n\in N$, we have $B'\res{S\times\{n\}}=0$ and
$B'\res{T_{(0,0,n)}M}=B_n'$. In fact, we will define $B'$ by an
explicit formula.

\sbr

Let us choose a coordinate system $\{x_i\}$ on $S$ centered at
$(0,0)$ (one can take $\{x_i\}=\{p,q\}$, but this is not important
in this step), and a coordinate system $\{y_j\}$ on $N$ centered
at $n_0$, so that $\{x_i,y_j\}$ is a coordinate system on $M$
centered at $m_0$. We denote the corresponding coordinate vector
fields by $\fs_i=\dd/\dd x_i$, $\fn_j=\dd/\dd y_j$. We then define
$B'$ by the formulas
\[
B'(\fs_i,\fs_k)=0; \qquad B'(\fs_i,\fn_j)(a,b,n) = B_n\bigl(
(\fs_i)_{(0,0,n)},(\fn_j)_{(0,0,n)} \bigr)
\]
(in particular, note that $B'(\fs_i,\fn_j)$ does not depend on the
coordinates $x_k$);
\[
B'(\fn_j,\fn_l) = \sum_i \Bigl[x_i\cdot \bigl(\fn_j
B'(\fs_i,\fn_l)-\fn_l B'(\fs_i,\fn_j) \bigr)\Bigr].
\]
By construction, $B'$ satisfies all the required pointwise
conditions, so we only have to check that $B'$ is closed. Since
$dB'$ is a differential $3$-form on $M$, it suffices to prove the
following identities:
\begin{equation}\label{e:1}
0=dB'(\fs_i,\fs_k,\fs_s)\overset{\text{def}}{=}\fs_i
B'(\fs_k,\fs_s) - \fs_k B'(\fs_i,\fs_s) + \fs_s B'(\fs_i,\fs_k);
\end{equation}
\begin{equation}\label{e:2}
0=dB'(\fs_i,\fs_k,\fn_j)\overset{\text{def}}{=}\fs_i
B'(\fs_k,\fn_j) - \fs_k B'(\fs_i,\fn_j) + \fn_j B'(\fs_i,\fs_k);
\end{equation}
\begin{equation}\label{e:3}
0=dB'(\fs_i,\fn_j,\fn_l)\overset{\text{def}}{=}\fs_i
B'(\fn_j,\fn_l) - \fn_j B'(\fs_i,\fn_l) + \fn_l B'(\fs_i,\fn_j);
\end{equation}
\begin{equation}\label{e:4}
0=dB'(\fn_j,\fn_l,\fn_t)\overset{\text{def}}{=}\fn_j
B'(\fn_l,\fn_t) - \fn_l B'(\fn_j,\fn_t) + \fn_t B'(\fn_j,\fn_l).
\end{equation}
The first two follow automatically from the definitions. The third
one follows from the definition of $B'(\fn_j,\fn_l)$ and the fact
that $\fs_i(x_k)=\de_{ik}$. Finally, the fourth identity follows
from a straightforward computation by substituting the definitions
of $B'(\fn_l,\fn_t)$, $B'(\fn_j,\fn_t)$ and $B'(\fn_j,\fn_l)$ into
\eqref{e:4} and using the fact that the vector fields $\fn_j$
commute with each other and annihilate $x_i$.

\subsection{Step 3}\label{ss:step3} We now complete the proof
outlined in \S\ref{ss:strategy}. Let us begin by exploring the
consequence of the Fundamental Theorem of Calculus in the context
of Lie derivatives. With the notation of \S\ref{ss:step1}, let
$\tau$ be a differential form on $M$ of arbitrary degree.
\begin{lem}\label{l:ftc}
For all $(a,b,n)\in S\times N=M$, we have
\begin{equation}\label{e:xp}
\tau_{(a,b,n)}=\phi_b^*\tau_{(a,0,n)} - \int_0^b \bigl(
\phi_{b-s}^*(\cL_{X_p}\tau)\bigr)_{(a,b,n)} ds
\end{equation}
and
\begin{equation}\label{e:xq}
\tau_{(a,b,n)} = \psi_{-a}^*\tau_{(0,b,n)} + \int_0^a
\bigl(\psi_{t-a}^*(\cL_{X_q}\tau)\bigr)_{(a,b,n)} dt.
\end{equation}
\end{lem}
The proof of this lemma is straightforward from the definition of
Lie derivative and the Fundamental Theorem of Calculus. Combining
\eqref{e:xp} and \eqref{e:xq}, we deduce that
\begin{eqnarray}\label{e:tau1}
\tau_{(a,b,n)} &=& \phi_b^*\psi_{-a}^*(\tau_{(0,0,n)}) + \phi_b^*
\int_0^a \bigl(\psi_{t-a}^*(\cL_{X_q}\tau)\bigr)_{(a,0,n)} dt -
\int_0^b \bigl(
\phi_{b-s}^*(\cL_{X_p}\tau)\bigr)_{(a,b,n)} ds \\
\label{e:tau2} &=& \psi_{-a}^*\phi_b^*(\tau_{(0,0,n)}) -
\psi_{-a}^* \int_0^b \bigl(
\phi_{b-s}^*(\cL_{X_p}\tau)\bigr)_{(0,b,n)} ds + \int_0^a
\bigl(\psi_{t-a}^*(\cL_{X_q}\tau)\bigr)_{(a,b,n)} dt.
\end{eqnarray}
In particular, Proposition \ref{p:xi}(3) now implies that
\begin{equation}\label{e:epsilon}
\eps_{(a,b,n)} = \phi_b^*\psi_{-a}^*(\eps_{(0,0,n)}) + \left\{
\phi_b^* \int_0^a \bigl(\psi_{t-a}^*(d\xi_q)\bigr)_{(a,0,n)} dt -
\int_0^b \bigl( \phi_{b-s}^*(d\xi_p)\bigr)_{(a,b,n)} ds \right\}
\Biggl\lvert_{E_{(a,b,n)}}.
\end{equation}
We now note that, due to the preparations of \S\S\ref{ss:step1}
and \ref{ss:step2}, the GC structure on $S\times N$ defined as the
product of the symplectic structure on $S$ and the induced GC
structure on $N$ corresponds to the $2$-form $\eps'$ on $E$
defined by
\[
\eps'_{(a,b,n)} = \phi_b^*\psi_{-a}^*\bigl( \eps_{(0,0,n)} \bigr).
\]
The proof will therefore be complete if we show that the (real)
$2$-form $B$ on $S\times N$ defined by
\[
B_{(a,b,n)} = \phi_b^* \int_0^a
\bigl(\psi_{t-a}^*(d\xi_q)\bigr)_{(a,0,n)} dt - \int_0^b \bigl(
\phi_{b-s}^*(d\xi_p)\bigr)_{(a,b,n)} ds
\]
is closed.

\sbr

Recall first from Proposition \ref{p:xi} that
$\xi_{\{f,g\}}=\cL_{X_f}(\xi_g)-\iota_{X_g}(d\xi_f)$ for all
$C^\infty$ functions $f,g$ on $M$; on the other hand, the
definition of the map $f\mapsto\xi_f$ implies that if $\{f,g\}$ is
a constant function on $M$, then $\xi_{\{f,g\}}=0$. We deduce that
\begin{align}
\label{e:xi1} \cL_{X_p}(\xi_q)&=\iota_{X_q}(d\xi_p), &
\cL_{X_q}(\xi_p)&=\iota_{X_p}(d\xi_q), \\
\label{e:xi2} \cL_{X_p}(\xi_p)&=\iota_{X_p}(d\xi_p), &
\cL_{X_q}(\xi_q)&=\iota_{X_q}(d\xi_q).
\end{align}

We now compute
\[
(\cL_{X_p}B)_{(a,b,n)} = \lim\limits_{\ga\to 0} \frac{1}{\ga}\cdot
\Bigl[ \phi_\ga^*\bigl(B_{(a,b-\ga,n)}\bigr) - B_{(a,b,n)} \Bigr],
\]
and since the flows $\phi_s$ and $\psi_t$ commute, we obtain
\[
\phi_\ga^*\bigl(B_{(a,b-\ga,n)}\bigr) = \phi_b^* \int_0^a
\bigl(\psi_{t-a}^*(d\xi_q)\bigr)_{(a,0,n)} dt - \int_0^{b-\ga}
\bigl( \phi_{b-s}^*(d\xi_p)\bigr)_{(a,b,n)} ds,
\]
whence
\[
\frac{1}{\ga}\cdot \Bigl[ \phi_\ga^*\bigl(B_{(a,b-\ga,n)}\bigr) -
B_{(a,b,n)} \Bigr] = \frac{1}{\ga}\cdot \int_{b-\ga}^b
\phi_{b-s}^*\bigl( (d\xi_p)_{(a,s,n)} \bigr) ds.
\]
The limit of this expression as $\ga\to 0$ is equal to
$(d\xi_p)_{(a,b,n)}$. Thus
\begin{equation}\label{e:lxp}
\cL_{X_p} B = d\xi_p.
\end{equation}
Similarly,
\[
(\cL_{X_q}B)_{(a,b,n)} = \lim\limits_{\ga\to 0} \frac{1}{\ga}\cdot
\Bigl[ \psi_\ga^*\bigl(B_{(a+\ga,b,n)}\bigr) - B_{(a,b,n)} \Bigr],
\]
and
\[
\psi_\ga^* B_{(a+\ga,b,n)} = \phi_b^* \int_0^{a+\ga}
\bigl(\psi_{t-a}^*(d\xi_q)\bigr)_{(a,0,n)} dt - \int_0^b
\psi_\ga^*\phi_{b-s}^*\bigl((d\xi_p)_{(a+\ga,s,n)}\bigr) ds,
\]
which leads to
\begin{equation}\label{e:intermediate}
(\cL_{X_q}B)_{(a,b,n)} = \phi_b^*(d\xi_q)_{(a,0,n)} - \int_0^b
\phi_{b-s}^*\Bigl( \bigl( \cL_{X_q}(d\xi_p) \bigr)_{(a,s,n)}
\Bigr) ds.
\end{equation}
However, we have, from \eqref{e:xi1} and Cartan's formula for
$\cL_{X_q}$,
\[
\cL_{X_q}(d\xi_p) = d\iota_{X_q}(d\xi_p) = d\cL_{X_p}(\xi_q) =
\cL_{X_p}(d\xi_q).
\]
Substituting this into \eqref{e:intermediate} and combining with
Lemma \ref{l:ftc}, we obtain
\begin{equation}\label{e:lxq}
\cL_{X_q} B = d\xi_q.
\end{equation}

\mbr

We now compute $\iota_{X_p}B$. We use the fact that contraction
commutes with integration of differential forms, and also that the
vector field $X_p$ is invariant under the flows $\phi_s$ and
$\psi_t$:
\begin{eqnarray*}
(\iota_{X_p}B)_{(a,b,n)} &=& \phi_b^* \int_0^a
\bigl(\psi_{t-a}^*(\iota_{X_p}d\xi_q)\bigr)_{(a,0,n)} dt -
\int_0^b \bigl(
\phi_{b-s}^*(\iota_{X_p}d\xi_p)\bigr)_{(a,b,n)} ds \\
&=& \phi_b^* \int_0^a
\bigl(\psi_{t-a}^*(\cL_{X_q}\xi_p)\bigr)_{(a,0,n)} dt - \int_0^b
\bigl(
\phi_{b-s}^*(\cL_{X_p}\xi_p)\bigr)_{(a,b,n)} ds \\
&=& \phi_b^* \Bigl( (\xi_p)_{(a,0,n)} - \psi_{-a}^*
\bigl((\xi_p)_{(0,0,n)}\bigr) \Bigr) + (\xi_p)_{(a,b,n)} -\phi_b^*
\bigl( (\xi_p)_{(a,0,n)} \bigr) \\
&=& (\xi_p)_{(a,b,n)} - \phi_b^*\psi_{-a}^*
\bigl((\xi_p)_{(0,0,n)}\bigr),
\end{eqnarray*}
where we have used \eqref{e:xi1}, \eqref{e:xi2} in the second
equality and Lemma \ref{l:ftc} in the third equality. However,
$(\xi_p)_{(0,0,n)}=0$. This follows from the fact that
$(\xi_p)_{(0,0,n)}$ depends only on the value of $dp$ at the point
$(0,0,n)$ and on the induced constant GC structure on
$T_{(0,0,n)}M$; on the other hand, after the preparations of
\S\S\ref{ss:step1}, \ref{ss:step2}, the constant GC structure on
$T_{(0,0,n)}M$ is the direct sum of the induced GC structure on
$T_n N$ and the symplectic GC structure on
$T_{(0,0,n)}\bigl(S\times\{n\}\bigr)$. Therefore
\begin{equation}\label{e:ixp}
\iota_{X_p} B = \xi_p.
\end{equation}

\sbr

Similarly,
\begin{eqnarray*}
(\iota_{X_q}B)_{(a,b,n)} &=& \phi_b^* \int_0^a
\bigl(\psi_{t-a}^*(\iota_{X_q}d\xi_q)\bigr)_{(a,0,n)} dt -
\int_0^b \bigl(
\phi_{b-s}^*(\iota_{X_q}d\xi_p)\bigr)_{(a,b,n)} ds \\
&=& \phi_b^* \int_0^a
\bigl(\psi_{t-a}^*(\cL_{X_q}\xi_q)\bigr)_{(a,0,n)} dt - \int_0^b
\bigl(
\phi_{b-s}^*(\cL_{X_p}\xi_q)\bigr)_{(a,b,n)} ds \\
&=& \phi_b^* \Bigl( (\xi_q)_{(a,0,n)} - \psi_{-a}^*
\bigl((\xi_q)_{(0,0,n)}\bigr) \Bigr) + (\xi_q)_{(a,b,n)} -\phi_b^*
\bigl( (\xi_q)_{(a,0,n)} \bigr) \\
&=& (\xi_q)_{(a,b,n)} - \phi_b^*\psi_{-a}^*
\bigl((\xi_q)_{(0,0,n)}\bigr) = (\xi_q)_{(a,b,n)},
\end{eqnarray*}
which forces
\begin{equation}\label{e:ixq}
\iota_{X_q} B = \xi_q.
\end{equation}

\sbr

Let us compare \eqref{e:lxp} and \eqref{e:ixp}. We can rewrite
\eqref{e:lxp} as $d(\iota_{X_p}B)+\iota_{X_p}(dB)=d\xi_p$, whence
\eqref{e:ixp} implies that $\iota_{X_p}(dB)=0$. Likewise,
\eqref{e:lxq} and \eqref{e:ixq} force $\iota_{X_q}(dB)=0$. But
$X_p,X_q$ span the tangent space to every horizontal leaf
$S\times\{n\}$ at every point. Hence, to show that $dB=0$, it
remains to check that the restriction of $dB$ to each transverse
slice $\{s\}\times N$ is equal to zero. By construction, the
restriction of $B$ itself to $\{(0,0\}\times N$ is zero. Let us
pick three arbitrary sections $Z_1,Z_2,Z_3$ of $\sN$ which commute
with $X_p$ and $X_q$. Then $(dB)(Z_1,Z_2,Z_3)=0$ along
$N=\{(0,0)\}\times N$, and furthermore
\[
\cL_{X_p}\bigl[ (dB)(Z_1,Z_2,Z_3) \bigr] =
(\cL_{X_p}(dB))(Z_1,Z_2,Z_3)=(d\cL_{X_p}B)(Z_1,Z_2,Z_3)=(dd\xi_p)(Z_1,Z_2,Z_3)=0,
\]
where we have used \eqref{e:lxp} and the fact that $X_p$ commutes
with each $Z_j$. Similarly, $\cL_{X_p}\bigl[ (dB)(Z_1,Z_2,Z_3)
\bigr]=0$. It follows that $(dB)(Z_1,Z_2,Z_3)=0$ everywhere on $M$
and completes the proof.

\subsection{Appendix: proof of Lemma \ref{l:subman}}\label{ss:prooflemma}
Let $(\sS,\om)$ denote the canonical symplectic foliation
associated to the GC structure on $M$, and recall from
\S\ref{ss:step1} that $\sS_0\subseteq\sS$ denotes the foliation
tangent to the leaves $S\times\{n\}$. Since $M=S\times N$ as
Poisson manifolds, it follows that at each point $(s,n)\in N$, the
tangent space $T_{(s,n)}
\bigl(S\times\{n\}\bigr)=\bigl(\sS_0\bigr)_{(s,n)}$ is orthogonal
to $T_{s,n}\bigl(\{s\}\times N\bigr)\cap\sS_{(s,n)}$ with respect
to $\om$. In particular, by Theorem \ref{t:gclinalg}(e), the
transverse slice $N$ satisfies the pointwise condition for being a
GC submanifold of $M$, and hence we must only show that $L_N$ is a
subbundle of $\ttsc{N}$.

\sbr

Since $L_N$ is the image of $\tilde{L}_N=L\res{N}\cap \Bigl(
\tc{N}\oplus \bigl(\tcs{M}\res{N}\bigr) \Bigr)$ under the
projection map $\tc{N}\oplus\bigl(\tcs{M}\res{N}\bigr)\to
\ttsc{N}$, it suffices to show that $\tilde{L}_N$ is a subbundle
of $\tc{N}\oplus\bigl(\tcs{M}\res{N}\bigr)$. Further, since
$\tilde{L}_N$ is defined as the intersection of two subbundles of
$\bigl(\ttsc{M}\bigr)\res{N}$, it suffices to show that
$\tilde{L}_N$ has constant rank on $N$. Considering the projection
of $\tilde{L}_N$ onto $\tc{N}$, we obtain a short exact sequence
\[
0 \rar{} \bigl(L\cap\tcs{M}\bigr)\res{N} \rar{} \tilde{L}_N \rar{}
\bigl( E\res{N}\cap\tc{N} \bigr) \rar{} 0.
\]
Now $L\cap\tcs{M}=\Ann_{\tcs{M}}(E)$, so
\begin{eqnarray*}
\rk \tilde{L}_N &=& \rk\bigl(\Ann_{\tcs{M}}(E)\res{N}\bigr) + \rk
\bigl(E\res{N}\bigr) - \rk \Bigl( \bigl(E\res{N}\bigr)\Big/ \bigl(
E\res{N}\cap \tc{N} \bigr) \Bigr) \\
&=& \dim M - \rk \Bigl( \bigl(E\res{N}+\tc{N}\bigr) \Big/ \tc{N}
\Bigr) \\
&=& \dim M - \rk \bigl( \tc{M}\res{N}\big/ \tc{N} \bigr) = \dim N;
\end{eqnarray*}
we have used the fact that $E\res{N}+\tc{N}=\tc{M}\res{N}$, which
follows from $\sS\res{N}+TN=TM\res{N}$. Thus, in fact, not only is
the rank of $\tilde{L}_N$ constant, but the projection map
$\tilde{L}_N\to L_N$ is an isomorphism (since $L_N$ has constant
rank equal to $\dim N$).

\section{Linearization of generalized complex
structures}\label{s:linear}

In this section we consider a GC structure $\cJ$ on a manifold $M$
such that the associated Poisson tensor has rank zero at a certain
point $m\in M$. Our goal is to describe a ``first-order
approximation'' to the GC structure in a neighborhood of $m$. We
will use the notation $\fm_m^2\subset\fm_m\subset C^\infty(M)$ in
the same sense as in \S\ref{ss:linearGCS}. Also, for $f\in
C^\infty(M)$, we will use the notation $(X_f,\xi_f)$ as defined in
Section \ref{s:poisson}. Let us assume that $\cJ$ is given by the
matrix \eqref{e:jmatrix}. Thus, by assumption, $\pi_m:T_m^* M\to
T_m M$ is the zero map. Hence, if we consider the induced constant
GC structure $\cJ_m$ on $T_m M$, its matrix has the form
\[
\cJ_m=\matr{J_m}{0}{\sg_m}{-J_m^*}.
\]
It is proved, for instance, in \cite{oren}, that a constant GC
structure of this form is a $B$-field transform of a complex GC
structure on $T_m M$. If $B_m\in\bigwedge^2 T_m^* M$ is any
two-form which transforms $\cJ_m$ into a complex GC structure, we
can extend $B_m$ to a differential $2$-form $B$ on a neighborhood
of $m$ in $M$ which is constant in the appropriate local
coordinates, and hence, a fortiori, is closed. Applying the
transformation defined by $B$ to the structure $\cJ$ reduces us to
the situation where $\sg_m=0$.

\sbr

We now assume that $\sg_m=0$ and explain what we mean by the
first-order approximation to $\cJ$ at the point $m$, proving
Theorem \ref{t:approx} at the same time. Let $\fg=\fm_m/\fm_m^2$
be the real Lie algebra which encodes the first-order
approximation to $\pi$ at $m$, as defined in \S\ref{ss:linearGCS}.
Thus the Lie bracket on $\fg$ is induced by the Poisson bracket on
$C^\infty(M)$ defined by $\pi$. We can also think of $\pi$ as a
$C^\infty(M)$-linear map from $\Ga(M,T^*M)$ to
$\fm_m\cdot\Ga(M,TM)$, which induces an $\bR$-linear map
\[
\Ga(M,T^*M)\big/\fm_m\cdot\Ga(M,T^*M) \rar{} \fm_m\cdot\Ga(M,TM)
\big/ \fm_m^2\cdot\Ga(M,TM).
\]
This map also encodes the first-order approximation to $\pi$. It
is then natural to define the first-order approximation to $\cJ$
to be the $\bR$-linear automorphism of
\[
\Bigl(\fm_m\cdot\Ga(M,TM) \big/ \fm_m^2\cdot\Ga(M,TM)\Bigr) \oplus
\Bigl( \Ga(M,T^*M)\big/\fm_m\cdot\Ga(M,T^*M) \Bigr)
\]
induced by $\cJ$. Note, however, that the map
\[
\fm_m\cdot\Ga(M,TM) \big/ \fm_m^2\cdot\Ga(M,TM) \rar{}
\Ga(M,T^*M)\big/\fm_m\cdot\Ga(M,T^*M)
\]
induced by $\sg$ clearly vanishes; moreover, since $J$ and
$K=-J^*$ determine each other, we can concentrate our attention on
the map
\[
\Ga(M,T^*M)\big/\fm_m\cdot\Ga(M,T^*M) \rar{}
\Ga(M,T^*M)\big/\fm_m\cdot\Ga(M,T^*M)
\]
induced by $-J^*$. Now the de Rham differential $d$ induces an
$\bR$-linear isomorphism
\[
\fm_m/\fm_m^2 \rar{\simeq} \Ga(M,T^*M)\big/\fm_m\cdot\Ga(M,T^*M),
\]
and we have, by definition $\xi_f=-J^*(df)$ for any $f\in\fm_m$.
By transport of structure, $-J^*$ induces an $\bR$-linear
automorphism of the Lie algebra $\fg$ which we will denote by $A$;
by construction, $A^2=-1$. To obtain further information on $A$,
we will study it from the point of view of the map
$f\mapsto\xi_f$.

\sbr

Let $f,g\in\fm_m$. Part (2) of Proposition \ref{p:xi} yields
\[
\xi_{\{f,g\}} = \cL_{X_f}(\xi_g)-\iota_{X_g}(d\xi_f) = \bigl[
\iota_{X_f}(d\xi_g) - \iota_{X_g}(d\xi_f) \bigr] +
d\bigl(\iota_{X_f}(\xi_g) \bigr).
\]
But $X_f=\pi(df)$ and $X_g=\pi(dg)$, which implies that the first
term vanishes modulo $\fm_m$. On the other hand, $\xi_g\equiv
d(Ag)$ modulo $\fm_m$, whence $\iota_{X_f}(\xi_g)\equiv \{f,Ag\}$
modulo $\fm_m^2$. Thus we conclude that $A\{f,g\}\equiv\{f,Ag\}$
modulo $\fm_m^2$, which is precisely the condition for $A$ to make
$\fg$ a {\em complex} Lie algebra. This completes the proof of
Theorem \ref{t:approx}.

\end{document}